\documentclass[12pt]{amsart}
\usepackage{tipa}
\usepackage{amssymb}
\usepackage{graphicx}
\theoremstyle{definition}

\theoremstyle{remark}

\numberwithin{equation}{section}

\begin{document}
\title[]{{\Large \bf Asymptotic behaviour for a class of subharmonic Functions
in a Half Space  $^{\ast}$ } }
\author{   Pan Guoshuang$^{1,2}$ and Deng Guantie$^{1, \ast\ast}$}%
\address{$^{1}$ Sch. Math. Sci. \& Lab. Math. Com. Sys. \\
   Beijing Normal University\\
    100875  Beijing, The People's Republic of China
    }%
\address{$^{2}$ Department of Public Basic Courses \\
   Beijing  Institute of Fashion and Technology \\
    100029  Beijing, The People's Republic of China
    }%
\email{denggt@bnu.edu.cn}%

\thanks{$\ast $ The project supported by NSFC (Grant No.10671022) and by RFDP (Grant
No.20060027023)
\endgraf $\ast\ast$:Corresponding author.}

\begin{abstract}
 \ \ {  A class of subharmonic  functions are proved to have the  growth estimates
 $u(x)= o(x_n^{1-\frac{\alpha}{p}}|x|^{\frac{\gamma}{p}+\frac{n-1}{q}-n+\frac{\alpha}{p}})$
 at infinity in the upper
 half space of ${\bf R}^{n}$, which generalizes the growth properties
 of analytic functions and harmonic functions. }
 \end{abstract}
\maketitle


 \section*{ 1. Introduction and Main Theorem}

\vspace{0.3cm}

  { Let ${\bf R}^{n} (n\geq3)$  denote the  $n$-dimensional Euclidean
space with points $x=(x_1,x_2,\cdots,x_{n-1},x_{n})=(x',x_n)$, where
$x' \in {\bf R}^{n-1}$ and $x_{n} \in {\bf R}$.  The boundary and
closure of an open  $\Omega$ of ${\bf R}^{n}$ are denoted by
$\partial{\Omega}$
 and $\overline{\Omega}$ respectively.
 The upper half-space $H$ is the set
 $H=\{x=(x',x_n)\in {\bf R}^{n}:\; x_n>0\}$, whose boundary is
 $\partial{H}$ .
    We  write $B(x,\rho)$ and $\partial B(x,\rho) $ for the open ball
    and the sphere of radius $\rho$  centered at $x$ in ${\bf R}^{n}$.
 We identify ${\bf R}^{n}$ with ${\bf R}^{n-1}\times {\bf R}$ and
${\bf R}^{n-1}$ with $ {\bf R}^{n-1}\times \{0\}$,
  with this convention we then have $ \partial {H}={\bf R}^{n-1}$,
  writing typical points $x,\ y \in {\bf R}^{n}$ as $x=(x',x_n),\
y=(y',y_n),$  where $x'=(x_1,x_2,\cdots,x_{n-1}),\
y'=(y_1,y_2,\cdots y_{n-1}) \in {\bf R}^{n-1}$ and putting
$$
x\cdot y=\sum_{j=1}^{n}x_jy_j=x'\cdot y'+x_ny_n,\ \ |x|=\sqrt{x\cdot
x},\ \ |x'|=\sqrt{x'\cdot x'}.
$$

  For $x\in{\bf R}^{n}\backslash\{0\}$, let([\textbf{10}])
$$
 E(x)=-r_n|x|^{2-n},
$$
where $|x|$ is the Euclidean norm, $r_n=\frac{1}{(n-2)\omega_{n}}$
and $\omega_{n}=\frac{2\pi^{\frac{n}{2}}}{\Gamma(\frac{n}{2})}$ is
the surface area of the unit sphere in ${\bf R}^{n} $. We know that
$E$ is
locally integrable in ${\bf R}^{n} $.\\
  The Green function $G(x,y)$ for the upper half space
 $H$ is given by([\textbf{10}])
$$
 G(x,y)=E(x-y)-E(x-y^{\ast}) \qquad x,y\in\overline{H} ,\  x\neq y,
 \eqno{(1.1)}
$$
where $^{\ast}$ denotes reflection in the boundary plane $\partial
H$ just as $y^{\ast}=(y_1,y_2,\cdots,y_{n-1},-y_n)$, then we define
the Poisson kernel $P(x,y')$ when $x\in H$ and $y'\in \partial H $
by
$$
 P(x,y')=-\frac{\partial G(x,y)}{\partial
 y_n}\bigg|_{y_n=0}=\frac{2x_n}{\omega_n|x-(y',0)|^n}.  \eqno{(1.2)}
$$

  The  Dirichlet problem of upper half space is to find a function
 $u$ satisfying
$$
 u\in C^2(H), \eqno{(1.3)}
$$
$$
 \Delta u=0,   x\in H, \eqno{(1.4)}
$$
$$
 \lim_{x\rightarrow x'}u(x)=f(x')\ {\rm nontangentially  \  a.e.}x'\in \partial H, \eqno{(1.5)}
$$
where $f$ is a measurable function of ${\bf R}^{n-1} $. The Poisson
integral of the upper half space is defined by
$$
u(x)=P[f](x)=\int_{{\bf R}^{n-1}}P(x,y')f(y')dy'.\eqno{(1.6)}
$$
 As we all know, the Poisson integral $P[f]$ exists if
$$
\int_{{\bf R}^{n-1}}\frac{|f(y')|}{1+|y'|^n} dy'<\infty.
$$
(see [\textbf{1,2}] and [\textbf{11}])In this paper, we will
consider measurable functions $f$ in ${\bf R}^{n-1}$ satisfying
$$
\int_{{\bf R}^{n-1}}\frac{|f(y')|^p}{(1+|y'|)^{\gamma}}
dy'<\infty.\eqno{(1.7)}
$$

   Siegel-Talvila([\textbf{5}]) have proved the
  following result:

\vspace{0.2cm}
 \noindent
{\bf Theorem A } Let $f$ be a measurable function in ${\bf R}^{n-1}$
satisfying (1.7). Then the harmonic function $v(x)$ defined by (1.6)
satisfies (1.3), (1.4), (1.5) and
$$
v(x)= o(x_n^{1-n}|x|^{n+m}) \quad  {\rm as}  \ |x|\rightarrow\infty.
$$

  In order to describe the asymptotic behaviour of subharmonic functions
in half-spaces([\textbf{8,9}] and [\textbf{10}]),
 we establish the following theorems.

\vspace{0.2cm}
 \noindent
{\bf Theorem 1} Let $1\leq p<\infty,\  \frac{1}{p}+\frac{1}{q}=1$
and
$$
-(n-1)(p-1)<\gamma <(n-1)+p  \quad  {\rm in \  case}  \ p>1;
$$
$$
 0<\gamma \leq n  \quad  {\rm in \  case}  \ p=1.
$$
If $f$ is a measurable function in ${\bf R}^{n-1}$ satisfying (1.4)
and  $v(x)$ is the harmonic function defined by (1.8), then there
exists $x_j\in H,\ \rho_j>0,$ such that
$$
\sum
_{j=1}^{\infty}\frac{\rho_j^{pn-\alpha}}{|x_j|^{pn-\alpha}}<\infty
\eqno{(1.8)}
$$
holds and
$$
v(x)=
o(x_n^{1-\frac{\alpha}{p}}|x|^{\frac{\gamma}{p}+\frac{n-1}{q}-n+\frac{\alpha}{p}})
\quad  {\rm as}  \ |x|\rightarrow\infty   \eqno{(1.9)}
$$
holds in $H-G$. where $ G=\bigcup_{j=1}^\infty B(x_j,\rho_j)$ and
$0< \alpha\leq n$.

\vspace{0.2cm}
 \noindent
 {\bf Remark 1 } If $\alpha=n$, $p=1$ and $\gamma=n$, then (1.8) is a finite sum,
 the set $G$ is the union of finite balls, so (1.9) holds in $H$. This is just
 the case $m=0$ of the result of Siegel-Talvila.

\vspace{0.2cm}
 \noindent
 {\bf Remark 2 } When $\gamma=-(n-1)(p-1)$, $p>1$, we have
$$
v(x)=
o(x_n^{1-\frac{\alpha}{p}}(\log|x|)^{\frac{1}{q}}|x|^{\frac{\gamma}{p}+\frac{n-1}{q}-n+\frac{\alpha}{p}})
\quad  {\rm as}  \ |x|\rightarrow\infty   \eqno{(1.10)}
$$
holds in $H-G$.

  Next, we will generalize Theorem 1 to subharmonic functions.

\vspace{0.2cm}
 \noindent
{\bf Theorem 2 } Let $p$ and $\gamma$ be as in Theorem 1. If $f$ is
a measurable function in ${\bf R}^{n-1}$ satisfying (1.7) and $\mu$
is a positive Borel measure satisfying
$$
\int_H\frac{y_n^p}{(1+|y|)^\gamma} d\mu(y)<\infty   \eqno{(1.10)}
$$
and
$$
\int_H\frac{1}{(1+|y|)^{n-1}} d\mu(y)<\infty.
$$
Write the subharmonic function
$$
u(x)= v(x)+h(x), \quad x\in H
$$
where $v(x)$ is the harmonic function defined by (1.8), $h(x)$ is
defined by
$$
h(x)= \int_H G(x,y)d\mu(y)
$$
and $G(x,y)$ is defined by (1.1). Then there exists $x_j\in H,\
\rho_j>0,$ such that (1.8) holds and
$$
u(x)=
o(x_n^{1-\frac{\alpha}{p}}|x|^{\frac{\gamma}{p}+\frac{n-1}{q}-n+\frac{\alpha}{p}})
\quad  {\rm as}  \ |x|\rightarrow\infty
$$
holds in $H-G$. where $ G=\bigcup_{j=1}^\infty B(x_j,\rho_j)$ and
$0< \alpha<2$.

\vspace{0.2cm}
 \noindent
 {\bf Remark 3 } When $\gamma=-(n-1)(p-1)$, $p>1$, we have
$$
u(x)=
o(x_n^{1-\frac{\alpha}{p}}(\log|x|)^{\frac{1}{q}}|x|^{\frac{\gamma}{p}+\frac{n-1}{q}-n+\frac{\alpha}{p}})
\quad  {\rm as}  \ |x|\rightarrow\infty   \eqno{(1.10)}
$$
holds in $H-G$.

\vspace{0.4cm}

\section*{2.   Proof of Theorem }

\vspace{0.3cm}

Let $\mu$ be a positive Borel measure  in ${\bf R}^n,\ \beta\geq0$,
the maximal function $M(d\mu)(x)$ of order $\beta$ is defined by
$$
M(d\mu)(x)=\sup_{ 0<r<\infty}\frac{\mu(B(x,r))}{r^\beta},
$$
then the maximal function $M(d\mu)(x):{\bf R}^n \rightarrow
[0,\infty)$ is lower semicontinuous, hence measurable. To see this,
for any $ \lambda >0 $, let $D(\lambda)=\{x\in{\bf
R}^{n}:M(d\mu)(x)>\lambda\}$. Fix $x \in D(\lambda)$, then there
exists
 $r>0$ such that $\mu(B(x,r))>tr^\beta$ for some $t>\lambda$, and
there exists $ \delta>0$ satisfying
$(r+\delta)^\beta<\frac{tr^\beta}{\lambda}$. If $|y-x|<\delta$, then
$B(y,r+\delta)\supset B(x,r)$, therefore $\mu(B(y,r+\delta))\geq
tr^\beta >\lambda(r+\delta)^\beta$. Thus $B(x,\delta)\subset
D(\lambda)$. This proves that $D(\lambda)$ is open for each
$\lambda>0$.

 In order to obtain the results, we
need these lemmas below:

\vspace{0.2cm}
 \noindent
{\bf Lemma 1 } Let $\mu$ be a positive Borel measure  in ${\bf
R}^n,\ \beta\geq0,\ \mu({\bf R}^n)<\infty,$ for any $ \lambda \geq
5^{\beta} \mu({\bf R}^n)$, set
$$
E(\lambda)=\{x\in{\bf R}^{n}:|x|\geq2,M(d\mu)(x) >
\frac{\lambda}{|x|^{\beta}}\}
$$
then  there exists $ x_j\in E(\lambda)\  ,\ \rho_j> 0,\
j=1,2,\cdots$, such that
$$
E(\lambda) \subset \bigcup_{j=1}^\infty B(x_j,\rho_j) \eqno{(2.1)}
$$
and
$$
\sum _{j=1}^{\infty}\frac{\rho_j^{\beta}}{|x_j|^{\beta}}\leq
\frac{3\mu({\bf R}^n)5^{\beta}}{\lambda} .\eqno{(2.2)}
$$
Proof: Let $E_k(\lambda)=\{x\in E(\lambda):2^k\leq |x|<2^{k+1}\}$,
then  for any $ x \in E_k(\lambda),$ there exists $ r(x)>0$, such
that $\mu(B(x,r(x))) >\lambda(\frac{r(x)}{|x|})^{\beta} $, therefore
$r(x)\leq 2^{k-1}$.
 Since $E_k(\lambda)$ can be covered
by
 the union of a family of balls $\{B(x,r(x)):x\in E_k(\lambda) \}$,
 by the Vitali Lemma([\textbf{6}]), there exists $  \Lambda_k\subset E_k(\lambda)$,
$\Lambda_k$ is at most countable, such that $\{B(x,r(x)):x\in
\Lambda_k \}$ are disjoint and
$$
E_k(\lambda) \subset
 \cup_{x\in \Lambda_k} B(x,5r(x)),
$$
so
$$
E(\lambda)=\cup_{k=1}^\infty E_k(\lambda) \subset \cup_{k=1}^\infty
\cup_{x\in \Lambda_k} B(x,5r(x)).\eqno{(2.3)}
$$

  On the other hand, note that $ \cup_{x\in \Lambda_k} B(x,r(x)) \subset \{x:2^{k-1}\leq
|x|<2^{k+2}\} $, so that
$$
 \sum_{x \in \Lambda_k}\frac{(5r(x))^{\beta}}{|x|^{\beta}}
\leq 5^\beta\sum_{x\in\Lambda_k}\frac{\mu(B(x,r(x)))}{\lambda} \leq
\frac{5^\beta}{\lambda} \mu\{x:2^{k-1}\leq |x|<2^{k+2}\}.
$$
Hence we obtain
$$
 \sum _{k=1}^{\infty}\sum
_{x \in \Lambda_k}\frac{(5r(x))^{\beta}}{|x|^{\beta}}
 \leq
 \sum _{k=1}^{\infty}\frac{5^\beta}{\lambda} \mu\{x:2^{k-1}\leq |x|<2^{k+2}\}
 \leq
\frac{3\mu({\bf R}^n)5^{\beta}}{\lambda}.
$$
  Rearrange $ \{x:x \in \Lambda_k,k=1,2,\cdots\} $ and $
\{5r(x):x \in \Lambda_k,k=1,2,\cdots\}
 $, we get $\{x_j\}$
and $\{\rho_j\}$ such that
 (2.1) and
(2.2) hold.

\vspace{0.2cm}
 \noindent
{\bf Lemma 2 }  The kernel $\frac{1}{|x-y|^n}$ has the following
estimates:\\
(1) If $|y|\leq \frac{|x|}{2}$, then $\frac{1}{|x-y|^n}\leq
\frac{2^n}{|x|^n}$;\\
(2) If $|y|> \frac{|x|}{2}$, then $\frac{1}{|x-y|^n}\leq
\frac{2^n}{|y|^n}$.\\

  Throughout the proof, $A$
denote various positive constants.

 \emph{Proof of Theorem 1}

 We prove only the case $p>1$; the proof of the case $p=1$ is similar. Suppose
\begin{eqnarray*}
&G_1& =\{y'\in {\bf R}^{n-1}: 1<|y'|\leq \frac{|x|}{2}\},\\
&G_2& =\{y'\in {\bf R}^{n-1}: \frac{|x|}{2}<|y'| \leq 2|x|\}, \\
&G_3& =\{y'\in {\bf R}^{n-1}: |y'|>2|x|\}, \\
&G_4& =\{y'\in {\bf R}^{n-1}: |y'|\leq 1\}. \\
\end{eqnarray*}

Define the measure $dm(y')$ by
$$
dm(y')=\frac{|f(y')|^p}{(1+|y'|)^{\gamma}} dy'
$$
  For any $\varepsilon >0$, there exists $R_\varepsilon >2$, such that
$$
\int_{|y'|\geq
R_\varepsilon}dm(y')\leq\frac{\varepsilon^p}{5^{pn-\alpha}}.
$$
For every Lebesgue measurable set $E \subset {\bf R}^{n-1}$  , the
measure $m^{(\varepsilon)}$ defined by $m^{(\varepsilon)}(E)
=m(E\cap\{x'\in{\bf R}^{n-1}:|x'|\geq R_\varepsilon\}) $ satisfies
$m^{(\varepsilon)}({\bf
R}^{n-1})\leq\frac{\varepsilon^p}{5^{pn-\alpha}}$, write
\begin{eqnarray*}
&v_1(x)& =\int_{G_1} P(x,y')f(y')
dy',\\
&v_2(x)&=\int_{G_2} P(x,y')f(y')
dy', \\
&v_3(x)&=\int_{G_3} P(x,y')f(y')
dy', \\
&v_4(x)&=\int_{G_4} P(x,y')f(y')
dy', \\
\end{eqnarray*}
then
$$
v(x) =v_1(x)+v_2(x)+v_3(x)+v_4(x). \eqno{(2.3)}
$$
Let $ E_1(\lambda)=\{x\in{\bf R}^{n}:|x|\geq2,\exists
t>0,m^{(\varepsilon)}(B(x,t)\cap{\bf R}^{n-1}
)>\lambda^p(\frac{t}{|x|})^{pn-\alpha}\}$, therefore, if $ |x|\geq
2R_\varepsilon$ and $x \notin E_1(\lambda)
 $, then we have
$$
\forall t>0,\ m^{(\varepsilon)}(B(x,t)\cap{\bf R}^{n-1}
)\leq\lambda^p(\frac{t}{|x|})^{pn-\alpha}.
$$
  First, if $\gamma >-(n-1)(p-1)$, then $\frac{\gamma q}{p}
  +(n-1)>0$.
  For $r >1$, we have
$$
v_1(x)=\int_{G_1\cap B(0,r) } P(x,y')f(y') dy' +\int_{G_1-B(0,r) }
P(x,y')f(y') dy'=v_{11}(x)+v_{12}(x).
$$
If $|x|>2r$, then we obtain by Lemma 2 (1) and Holder's inequality
\begin{eqnarray*}
|v_{11}(x)|
&\leq& \int_{B(0,r)-B(0,1)}\frac{ 2x_n}{\omega_n}\frac{2^n}{|x|^n}|f(y')| dy' \\
&\leq& \frac{ 2^{n+1}}{\omega_n}\frac{x_n}{|x|^n}
\bigg(\int_{B(0,r)-B(0,1)}\frac{|f(y')|^p}{|y'|^\gamma}
dy'\bigg)^{1/p}\bigg(\int_{B(0,r)-B(0,1)}|y'|^{\frac{\gamma q}{p}}
dy'\bigg)^{1/q}, \\
\end{eqnarray*}
since
$$
\int_{B(0,r)-B(0,1)}|y'|^{\frac{\gamma q}{p}} dy'\leq
\omega_{n-1}\frac{1}{\frac{\gamma q}{p}+n-1}r^{\frac{\gamma
q}{p}+n-1} ,
$$
so that
$$
v_{11}(x)= o(x_n|x|^{\frac{\gamma}{p}+\frac{n-1}{q}-n}) \quad {\rm
as}  \ |x|\rightarrow\infty .  \eqno{(1.9)}
$$
Moreover, we have similarly

\begin{eqnarray*}
|v_{12}(x)|
&\leq& \frac{ 2^{n+1}}{\omega_n}\frac{x_n}{|x|^n}
\bigg(\int_{G_1-B(0,r)}\frac{|f(y')|^p}{|y'|^\gamma}
dy'\bigg)^{1/p}\bigg(\int_{G_1-B(0,r)}|y'|^{\frac{\gamma q}{p}}
dy'\bigg)^{1/q} \\
&\leq& Ax_n|x|^{\frac{\gamma}{p}+\frac{n-1}{q}-n}
\bigg(\int_{G_1-B(0,r)}\frac{|f(y')|^p}{|y'|^\gamma}
dy'\bigg)^{1/p}, \\
\end{eqnarray*}
which implies by artitrariness of $r$ that
$$
v_{12}(x)= o(x_n|x|^{\frac{\gamma}{p}+\frac{n-1}{q}-n}) \quad {\rm
as}  \ |x|\rightarrow\infty .  \eqno{(1.9)}
$$
If $\gamma >-(n-1)(p-1)$, then $\frac{\gamma q}{p}
  +(n-1)>0$, so that we obtain by Holder's inequality
\begin{eqnarray*}
|v_2(x)|
&\leq& \frac{2x_n}{\omega_n}
\bigg(\int_{G_2}\frac{|f(y')|^p}{|x-(y',0)|^{pn}|y'|^\gamma}
dy'\bigg)^{1/p}\bigg(\int_{G_2}|y'|^{\frac{\gamma q}{p}}
dy'\bigg)^{1/q} \\
&\leq& Ax_n|x|^{\frac{\gamma}{p}+\frac{n-1}{q}}
\bigg(\int_{G_2}\frac{|f(y')|^p}{|x-(y',0)|^{pn}|y'|^\gamma}
dy'\bigg)^{1/p},\\
\end{eqnarray*}
since
\begin{eqnarray*}
\int_{G_2}\frac{|f(y')|^p}{|x-(y',0)|^{pn}|y'|^\gamma} dy'
&\leq& \int_\frac{x_n}{2}^{3|x|}
\frac{2^\gamma+1}{t^{pn}} dm_x^{(\varepsilon)}(t) \\
&\leq& \frac{\varepsilon^p}{
|x|^{pn}}(2^\gamma+1)\bigg(\frac{1}{3^\alpha}+
\frac{pn}{\alpha}\bigg)\frac{|x|^\alpha}{x_n^\alpha}, \\
\end{eqnarray*}
where  $m_x^{(\varepsilon)}(t)=\int_{|y'-x| \leq t}
dm^{(\varepsilon)}(y')$.\\
Hence we have
$$
v_2(x)=
o(x_n^{1-\frac{\alpha}{p}}|x|^{\frac{\gamma}{p}+\frac{n-1}{q}-n+\frac{\alpha}{p}})
\quad  {\rm as}  \ |x|\rightarrow\infty.
$$
If $\gamma <(n-1)+p$, then $(\frac{\gamma}{p}-n)q
  +(n-1)<0$, so that we obtain by Lemma 2 (2) and Holder's inequality

\begin{eqnarray*}
|v_3(x)|
&\leq& \int_{G_3}\frac{ 2x_n}{\omega_n}\frac{2^n}{|y'|^n}|f(y')| dy' \\
&\leq& \frac{ 2^{n+1}}{\omega_n}x_n
\bigg(\int_{G_3}\frac{|f(y')|^p}{|y'|^\gamma}
dy'\bigg)^{1/p}\bigg(\int_{G_3}|y'|^{(\frac{\gamma }{p}-n)q}
dy'\bigg)^{1/q} \\
&\leq& Ax_n|x|^{\frac{\gamma}{p}+\frac{n-1}{q}-n}
\bigg(\int_{G_3}\frac{|f(y')|^p}{|y'|^\gamma}
dy'\bigg)^{1/p},\\
\end{eqnarray*}
so that
$$
v_3(x)= o(x_n|x|^{\frac{\gamma}{p}+\frac{n-1}{q}-n}) \quad {\rm as}
\ |x|\rightarrow\infty .  \eqno{(1.9)}
$$
Finally, by Lemma 2 (1), we obtain
$$
|v_4(x)|\leq \frac{ 2^{n+1}}{\omega_n}\frac{x_n}{|x|^n}
\int_{G_4}{|f(y')|}dy',
$$
so that we have by $\gamma >-(n-1)(p-1)$
$$
v_4(x)= o(x_n|x|^{\frac{\gamma}{p}+\frac{n-1}{q}-n}) \quad {\rm as}
\ |x|\rightarrow\infty .  \eqno{(1.9)}
$$
  Thus, by collecting (2.5), (2.6), (2.7), (2.8), (2.9), (2.10)and
(2.11), there exists a positive constant $A$ independent of
$\varepsilon$, such that if $ |x|\geq 2R_\varepsilon$ and $\  x
\notin E_1(\varepsilon)$, we have
$$
|v(x)|\leq A\varepsilon
x_n^{1-\frac{\alpha}{p}}|x|^{\frac{\gamma}{p}+\frac{n-1}{q}-n+\frac{\alpha}{p}}.
$$

 Let $\mu_\varepsilon$ be a measure in ${\bf R}^n$ defined by
$ \mu_\varepsilon(E)= m^{(\varepsilon)}(E\cap{\bf R}^{n-1})$ for
every measurable set $E$ in ${\bf R}^n$.Take
$\varepsilon=\varepsilon_p=\frac{1}{2^{p+2}}, p=1,2,3,\cdots$, then
there exists a sequence $ \{R_p\}$: $1=R_0<R_1<R_2<\cdots$ such that
$$
\mu_{\varepsilon_p}({\bf R}^n)=\int_{|y'|\geq
R_p}dm(y')<\frac{\varepsilon_p^p}{5^{pn-\alpha}}.
$$
Take $\lambda=3\cdot5^{pn-\alpha}\cdot2^p\mu_{\varepsilon_p}({\bf
R}^n)$ in Lemma 1, then there exists $x_{j,p}$ and $ \rho_{j,p}$,
where $R_{p-1}\leq |x_{j,p}|<R_p$, such that
$$
\sum _{j=1}^{\infty}(\frac{\rho_{j,p}}{|x_{j,p}|})^{pn-\alpha} \leq
\frac{1}{2^{p}}.
$$
if $R_{p-1}\leq |x|<R_p$ and $x\notin G_p=\cup_{j=1}^\infty
B(x_{j,p},\rho_{j,p})$, we have
$$
|v(x)|\leq
A\varepsilon_px_n^{1-\frac{\alpha}{p}}|x|^{\frac{\gamma}{p}+\frac{n-1}{q}-n+\frac{\alpha}{p}},
$$
Thereby
$$
\sum _{p=1}^{\infty}
\sum_{j=1}^{\infty}(\frac{\rho_{j,p}}{|x_{j,p}|})^{pn-\alpha} \leq
\sum _{p=1}^{\infty}\frac{1}{2^{p}}=1<\infty.
$$

  Set $ G=\cup_{p=1}^\infty G_p$, thus Theorem 1 holds.

 \emph{Proof of Theorem 2}

 We prove only the case $p>1$; the remaining case $p=1$ can be proved similarly.
Suppose
\begin{eqnarray*}
&F_1& =\{y\in H: 1<|y|\leq \frac{|x|}{2}\},\\
&F_2& =\{y\in H: \frac{|x|}{2}<|y| \leq 2|x|\}, \\
&F_3& =\{y\in H: |y|>2|x|\}, \\
&F_4& =\{y\in H: |y|\leq 1\}. \\
\end{eqnarray*}

  Define the measure $dn(y)$ by
$$
dn(y)=\frac{y_n^p}{(1+|y'|)^{\gamma}} d\mu(y)
$$

  For any $\varepsilon >0$, there exists $R_\varepsilon >2$, such that
$$
\int_{|y|\geq
R_\varepsilon}dn(y)<\frac{\varepsilon^p}{5^{pn-\alpha}}.
$$
For every Lebesgue measurable set $E \subset {\bf R}^{n}$, the
measure $n^{(\varepsilon)}$ defined by $n^{(\varepsilon)}(E)
=n(E\cap\{y\in H:|y|\geq R_\varepsilon\}) $ satisfies
$n^{(\varepsilon)}(H)\leq\frac{\varepsilon^p}{5^{pn-\alpha}}$, write
\begin{eqnarray*}
&h_1(x)& = \int_{F_1} G(x,y)d\mu(y),\\
&h_2(x)&=\int_{F_2} G(x,y)d\mu(y),\\
&h_3(x)&=\int_{F_3} G(x,y)d\mu(y), \\
&h_4(x)&=\int_{F_4} G(x,y)d\mu(y) \\
\end{eqnarray*}
then
$$
h(x) =h_1(x)+h_2(x)+h_3(x)+h_4(x). \eqno{(2.10)}
$$
Let $ E_2(\lambda)=\{x\in{\bf R}^{n}:|x|\geq2,\exists
t>0,n^{(\varepsilon)}(B(x,t)\cap H
)>\lambda^p(\frac{t}{|x|})^{pn-\alpha}\}, $ therefore, if $ |x|\geq
2R_\varepsilon$ and $x\notin E_2(\lambda)
 $, then we have
$$
\forall t>0, \ n^{(\varepsilon)}(B(x,t)\cap H
)\leq\lambda^p(\frac{t}{|x|})^{pn-\alpha}.
$$
First, note that
$$
|G(x,y)|=|E(x-y)-E(x-y^{\ast})|\leq \frac{2x_ny_n}{\omega_n|x-y|^n}.
\eqno{(2.11)}
$$
If $\gamma >-(n-1)(p-1)$, then $\frac{\gamma q}{p}
  +(n-1)>0$.
  For $r >1$, we have
$$
h_1(x)=\int_{F_1\cap B(0,r) } -G(x,y) d\mu(y) +\int_{F_1-B(0,r) }
 -G(x,y) d\mu(y)=h_{11}(x)+h_{12}(x)
$$
If $|x|>2r$, then we obtain by Lemma 2 (1), (2.11) and Holder's
inequality
\begin{eqnarray*}
|h_{11}(x)|
&\leq& \int_{B(0,r)-B(0,1)}\frac{2x_ny_n}{\omega_n|x-y|^n} d\mu(y) \\
&\leq& \int_{B(0,r)-B(0,1)}\frac{2x_ny_n}{\omega_n}\frac{2^n}{|x|^n} d\mu(y) \\
&\leq& \frac{ 2^{n+1}}{\omega_n}\frac{x_n}{|x|^n}
\bigg(\int_{B(0,r)-B(0,1)}\frac{y_n^p}{|y|^\gamma}
d\mu(y)\bigg)^{1/p}\bigg(\int_{B(0,r)-B(0,1)}|y|^{\frac{\gamma
q}{p}}
d\mu(y)\bigg)^{1/q}, \\
\end{eqnarray*}
since
$$
\int_{B(0,r)-B(0,1)}|y|^{\frac{\gamma q}{p}} d\mu(y)\leq
 2^{n-1}r^{\frac{\gamma
q}{p}+n-1}\int_{H}\frac{1}{(1+|y|)^{n-1}} d\mu(y),
$$
\begin{eqnarray*}
\int_{B(0,r)-B(0,1)}|y|^{\frac{\gamma q}{p}} d\mu(y)
&=&  \int_{B(0,r)-B(0,1)}|y|^{\frac{\gamma q}{p}+n-1}\frac{1}{|y|^{n-1}} d\mu(y) \\
&\leq& 2^{n-1}\int_{H}\frac{1}{(1+|y|)^{n-1}} d\mu(y)r^{\frac{\gamma
q}{p}+n-1} ,
\end{eqnarray*}
so that
$$
h_{11}(x)= o(x_n|x|^{\frac{\gamma}{p}+\frac{n-1}{q}-n}) \quad {\rm
as}  \ |x|\rightarrow\infty .  \eqno{(1.9)}
$$
Moreover, we have similarly
\begin{eqnarray*}
|h_{12}(x)|
&\leq& \frac{ 2^{n+1}}{\omega_n}\frac{x_n}{|x|^n}
\bigg(\int_{F_1-B(0,r)}\frac{y_n^p}{|y|^\gamma}
d\mu(y)\bigg)^{1/p}\bigg(\int_{F_1-B(0,r)}|y|^{\frac{\gamma q}{p}}
d\mu(y)\bigg)^{1/q} \\
&\leq& Ax_n|x|^{\frac{\gamma}{p}+\frac{n-1}{q}-n}
\bigg(\int_{F_1-B(0,r)}\frac{y_n^p}{|y|^\gamma} d\mu(y)\bigg)^{1/p},\\
\end{eqnarray*}
which implies by artitrariness of $r$ that
$$
h_{12}(x)= o(x_n|x|^{\frac{\gamma}{p}+\frac{n-1}{q}-n}) \quad {\rm
as} \ |x|\rightarrow\infty .  \eqno{(1.9)}
$$
If $\gamma >-(n-1)(p-1)$, then $\frac{\gamma q}{p}
  +(n-1)>0$, so that we obtain by Holder's inequality
\begin{eqnarray*}
|h_2(x)|
&\leq& \bigg(\int_{F_2}\frac{|G(x,y)|^p}{|y|^\gamma}
d\mu(y)\bigg)^{1/p}\bigg(\int_{F_2}|y|^{\frac{\gamma q}{p}}
d\mu(y)\bigg)^{1/q} \\
&\leq& \bigg(\int_{F_2}\frac{|G(x,y)|^p}{y_n^p}(2^\gamma +1)
dn(y)\bigg)^{1/p}\bigg(\int_{F_2}|y|^{\frac{\gamma q}{p}}
d\mu(y)\bigg)^{1/q}\\
&\leq& A|x|^{\frac{\gamma}{p}+\frac{n-1}{q}}\bigg(\int_{F_2}\frac{|G(x,y)|^p}{y_n^p}dn(y)\bigg)^{1/p},\\
\end{eqnarray*}
since
\begin{eqnarray*}
\int_{F_2}\frac{|G(x,y)|^p}{y_n^p} dn(y)
&\leq&  \int_{|y-x|\leq 3|x|}\frac{|G(x,y)|^p}{y_n^p} dn^{(\varepsilon)}(y) \\
&=&  \int_{|y-x|\leq \frac{x_n}{2}}\frac{|G(x,y)|^p}{y_n^p}
dn^{(\varepsilon)}(y)+
 \int_{\frac{x_n}{2}<|y-x|\leq 3|x|}\frac{|G(x,y)|^p}{y_n^p} dn^{(\varepsilon)}(y)\\
&=& h_{21}(x)+h_{22}(x),
\end{eqnarray*}
so that
\begin{eqnarray*}
h_{21}(x)
&\leq& \int_{|y-x|\leq
\frac{x_n}{2}}\bigg(\frac{2}{(n-2)\omega_nx_n|x-y|^{(n-2)}}\bigg)^p
dn^{(\varepsilon)}(y) \\
&=& \bigg(\frac{2}{(n-2)\omega_nx_n}\bigg)^p\int_0^{\frac{x_n}{2}}
\frac{1}{t^{p(n-2)}} dn_x^{(\varepsilon)}(t) \\
&\leq&
\bigg(\frac{2}{(n-2)\omega_n}\bigg)^p\frac{np-\alpha}{(2p-\alpha)2^{2p-\alpha}}
\varepsilon^p\frac{x_n^{p-\alpha}}{|x|^{np-\alpha}}.\\
\end{eqnarray*}
Moreover, we have by (2.11)
\begin{eqnarray*}
h_{22}(x)
&\leq& \int_{\frac{x_n}{2}<|y-x|\leq 3|x|
}\bigg(\frac{2x_n}{(n-2)\omega_n|x-y|^{n}}\bigg)^p
dn^{(\varepsilon)}(y) \\
&=& \bigg(\frac{2x_n}{\omega_n}\bigg)^p\int_{\frac{x_n}{2}}^{3|x|}
\frac{1}{t^{pn}} dn_x^{(\varepsilon)}(t) \\
&\leq&
\bigg(\frac{2}{\omega_n}\bigg)^p(\frac{1}{3^\alpha}+\frac{np2^\alpha}{\alpha})
\varepsilon^p\frac{x_n^{p-\alpha}}{|x|^{np-\alpha}},\\
\end{eqnarray*}
where  $n_x^{(\varepsilon)}(t)=\int_{|y-x| \leq t}
dn^{(\varepsilon)}(y)$.\\
Hence we have
$$
h_2(x)=
o(x_n^{1-\frac{\alpha}{p}}|x|^{\frac{\gamma}{p}+\frac{n-1}{q}-n+\frac{\alpha}{p}})
\quad  {\rm as}  \ |x|\rightarrow\infty.
$$
If $\gamma <(n-1)+p$, then $(\frac{\gamma}{p}-n)q
  +(n-1)<0$, so that we obtain by Lemma 2 (2), (2.11) and Holder's inequality

\begin{eqnarray*}
|h_3(x)|
&\leq& \int_{F_3}\frac{2x_ny_n}{\omega_n|x-y|^n} d\mu(y) \\
&\leq& \int_{F_3}\frac{ 2x_ny_n}{\omega_n}\frac{2^n}{|y|^n} d\mu(y) \\
&\leq& \frac{ 2^{n+1}}{\omega_n}x_n
\bigg(\int_{F_3}\frac{y_n^p}{|y|^\gamma}
d\mu(y)\bigg)^{1/p}\bigg(\int_{F_3}|y|^{(\frac{\gamma }{p}-n)q}
d\mu(y)\bigg)^{1/q} \\
&\leq& Ax_n|x|^{\frac{\gamma}{p}+\frac{n-1}{q}-n}
\bigg(\int_{F_3}\frac{y_n^p}{|y|^\gamma} d\mu(y)\bigg)^{1/p},\\
\end{eqnarray*}
so that
$$
h_3(x)= o(x_n|x|^{\frac{\gamma}{p}+\frac{n-1}{q}-n}) \quad {\rm as}
\ |x|\rightarrow\infty .  \eqno{(1.9)}
$$
Finally, by Lemma 2 (1) and (2.11), we obtain
$$
|h_4(x)|\leq \int_{F_4}\frac{2x_ny_n}{\omega_n|x-y|^n} d\mu(y) \leq
\frac{ 2^{n+1}}{\omega_n}\frac{x_n}{|x|^n} \int_{F_4}y_n d\mu(y),
$$
so that we have by $\gamma >-(n-1)(p-1)$
$$
h_4(x)= o(x_n|x|^{\frac{\gamma}{p}+\frac{n-1}{q}-n}) \quad {\rm as}
\ |x|\rightarrow\infty .  \eqno{(1.9)}
$$

  Thus, by collecting (2.12), (2.13), (2.15), (2.16),
(2.17), (2.18), (2.19) and (2.20), there exists a positive constant
$A$ independent of $\varepsilon$, such that if $ |x|\geq
2R_\varepsilon$ and $\  x \notin E_2(\varepsilon)$, we have
$$
 |h(x)|\leq A\varepsilon x_n^{1-\frac{\alpha}{p}}|x|^{\frac{\gamma}{p}+\frac{n-1}{q}-n+\frac{\alpha}{p}}.
$$

  Similarly, if $x\notin G$, we have
$$
h(x)=
o(x_n^{1-\frac{\alpha}{p}}|x|^{\frac{\gamma}{p}+\frac{n-1}{q}-n+\frac{\alpha}{p}})\quad
{\rm as} \ |x|\rightarrow\infty. \eqno{(2.21)}
$$

by (1.11) and  (2.21), we obtain
$$
u(x)=v(x)+h(x)=
o(x_n^{1-\frac{\alpha}{p}}|x|^{\frac{\gamma}{p}+\frac{n-1}{q}-n+\frac{\alpha}{p}})\quad
{\rm as} \ |x|\rightarrow\infty
$$
hold in $H-G$, thus we complete the proof of Theorem 2.

\begin{center}

\end{center}

\end{document}